\documentclass[12pt]{amsart}

\usepackage{amssymb}
\usepackage{latexsym}
\usepackage{verbatim}
\usepackage{fullpage}
\usepackage{pslatex}

\input {cyracc.def}
\font\ququ=cmr10 scaled \magstep1
 \font\tencyr=wncyr10 
\font\tencyi=wncyi10 
\font\tencysc=wncysc10 
\def\rus{\tencyr\cyracc}
\def\rusi{\tencyi\cyracc}
\def\rusc{\tencysc\cyracc}

\newcommand{\re}[1]{\textrm  (\ref{#1})}

\renewenvironment{proof}
{\noindent {\sl Proof.}\quad }{\hfill
$\square$ \vskip1.1ex\noindent }

\newenvironment{proof*}
{\noindent {\sl Proof.}\quad }{\hfill
$\square$}

\renewcommand{\theequation}{\thesection .\arabic{equation}}
\renewcommand{\thesubsubsection}{\theequation .\arabic{subsubsection}}

\catcode`\@=\active
\catcode`\@=11
\def\@eqnnum{\hbox to
.01pt{}\rlap{\hskip-\displaywidth(\mathbf{\theequation})}}
\catcode`\@=12

\newenvironment{s}[1]
{ \vskip1.2ex \refstepcounter{equation}
\noindent {\bf \theequation\enspace #1.} \begin{sl}}{\end{sl}
\vskip1.1ex\noindent }

\newenvironment{rem}[1]
{ \vskip1.2ex \refstepcounter{equation}
\noindent {\bf \theequation\enspace {#1}.} }{ \vskip1.1ex\noindent }

\newenvironment{subs}[1]
{\vskip1.2ex \refstepcounter{equation}
\noindent {\bf (\theequation)\quad #1.} }{\quad}

\newcommand {\be}{{\frak b}}
\newcommand {\ce}{{\frak c}}

\newcommand {\g}{{\frak g}}

\newcommand {\te}{{\frak t}}
\newcommand {\ut}{{\frak u}}

\newcommand {\gA}{{\goth A}}
\newcommand {\gC}{{\goth C}}
\newcommand {\gH}{{\goth H}}

\newcommand {\ap}{\alpha}

\newcommand {\lb}{\lambda}
\newcommand {\vp}{\varphi}

\newcommand {\HW}{\widehat W}
\newcommand {\HP}{\widehat\Pi}
\newcommand {\HD}{\widehat\Delta}

\newcommand {\ca}{{\mathcal A}}

\newcommand {\cF}{{\mathcal F}}

\newcommand {\VV}{{\Bbb V}}
\newcommand {\UU}{{\Bbb U}}

\newcommand {\asst}{{\ast\,}}

\newcommand {\GR}[2]{{\textrm{{\bf #1}}}_{#2}}

\newcommand {\ov}{\overline}
\newcommand {\un}{\underline}

\newcommand {\AD}{{\frak Ad}}
\newcommand {\AN}{{\frak An}}

\newcommand {\beq}{\begin{equation}}
\newcommand {\eeq}{\end{equation}}

\newcommand{\Cat}{\mathrm{Cat}}
\newcommand{\cats}{\mathrm{Cat}_s(\Delta)}
\newcommand{\cts}[1]{\mathrm{Cat}_s^{#1}(\Delta)}
\newcommand{\catsm}{\mathrm{Cat}_s^m(\Delta)}
\newcommand{\curge}{\succcurlyeq}
\newcommand{\curle}{\preccurlyeq}

\newcommand{\vts}{\VV(\theta_s)}

\newfam\Bbbfam\newfam\eufam\newfam\eusfam%
\font\Bbbfont=msbm10 scaled 1200%
\font\olala=msam10 scaled 1200%
\font\frak=eufm10 scaled 1400%
\font\Bbbsmallfont=msbm8%
\textfont\Bbbfam=\Bbbfont\scriptfont\Bbbfam=\Bbbsmallfont
\font\euzw=eufm10 scaled 1200%
\font\euac=eufm7 scaled 1200%
\font\euacc=eufm7 scaled 1000%
\textfont\eufam=\euzw\scriptfont\eufam=\euac%
\scriptscriptfont\eufam=\euacc%
\font\euszw=eusm10 scaled 1200%
\font\eusac=eusm7 scaled 1200%
\font\eusacc=eusm7 scaled 1000%
\textfont\eusfam=\euszw\scriptfont\eusfam=\eusac%
\scriptscriptfont\eusfam=\eusacc%
\def\frak{\fam\eufam}%
\def\goth{\fam\eusfam}%
\def\Bbb{\fam\Bbbfam}%

\def\varnothing{\hbox {\Bbbfont\char'077}}
\def\square{\hbox {\olala\char"03}}

\def\cyeq{\hbox {\olala\char'064}}

\begin{document}
\setlength{\parskip}{2pt plus 4pt minus 0pt}
\hfill {\scriptsize \today}
\vskip1ex
\vskip1ex

\title[Short antichains]{Short antichains in root
systems, semi-Catalan arrangements, and ${B}$-stable subspaces
}
\author{\sc Dmitri I. Panyushev}
\thanks{This research was supported in part by R.F.B.I. Grants no.
01--01--00756 and  02--01--01041}
\maketitle
\begin{center}
{\footnotesize
{\it Independent University of Moscow,
Bol'shoi Vlasevskii per. 11 \\
121002 Moscow, \quad Russia \\ e-mail}: {\tt panyush@mccme.ru }\\
}
\end{center}

\noindent

Let $G$ be a complex simple algebraic group with Lie algebra $\g$.
Fix a Borel subalgebra $\be$. An ideal of $\be$ is called
{\sf ad}-{\it nilpotent\/}, if it is contained in $[\be,\be]$.
The goal of this paper is to present a refinement of the
enumerative theory of
{\sf ad}-nilpotent ideals in the case, where $\g$ has roots of different
length.
\\[.4ex]
Let ${\AD}$  denote the set of all {\sf ad}-nilpotent ideals
of $\be$. Any $\ce\in\AD$ is completely determined by the corresponding
set of roots. The minimal roots in this set are called the
{\it generators\/} of an ideal. The collection of generators of an ideal
form an {\it antichain\/}
in the poset of positive roots, and the whole theory can be expressed
in the combinatorial language, in terms of antichains.
An antichain is called {\it strictly positive\/}, if it contains no simple
roots. Enumerative results for all and strictly positive antichains
were recently obtained in the work Athanasiadis, Cellini-Papi, Sommers,
and this author
\cite{ath02},\,\cite{ath03},\,\cite{CP1},\,\cite{CP2},\,\cite{duality},\,%
\cite{eric}.
\\[.4ex]
There are two different theoretical approaches to describing
(enumerating)
antichains. The first approach consists of constructing a bijection
between antichains and the coroot lattice points lying in
a certain simplex. An important intermediate step here is a
bijection between antichains and the so-called {\it minimal\/}
elements of the affine Weyl group, $\HW$.
It turns out that the simplex obtained is "equivalent" to a dilation of the
fundamental alcove of $\HW$, so that the problem
of counting the coroot lattice points in it can be resolved.
For strictly positive antichains, one constructs another bijection and
another simplex, and the respective elements of $\HW$ are called
{\it maximal\/}; yet, everything is quite similar.
The second approach uses the Shi bijection between the {\sf ad}-nilpotent
ideals (or antichains) and the dominant regions of the
Catalan arrangement. Under this bijection,
the strictly positive antichains correspond to the bounded regions.
There is a powerful result of
Zaslavsky allowing one to compute the number of all and bounded regions, if
the characteristic polynomial of the arrangement is known.
Since the characteristic polynomial of the Catalan arrangement
is recently computed in \cite{ath02}, the result follows.
\\[.4ex]
If $\g$ has roots of different length, one can distinguish the length
of elements occurring in antichains. We say that an antichain is
{\it short\/}, if it consists of only short roots.
This notion has a natural representation-theoretic incarnation:
the short antichains are in a one-to-one correspondence with
the $\be$-stable subspaces, without nonzero semisimple elements, in
the little adjoint $G$-module. A short analogue of strictly positive
antichains, {\it strictly $s$-positive\/} antichains, is also
defined. We are able to carry the above two
approaches over to the short antichains. First, we
introduce and characterize suitable elements of $\HW$
($s$-{\it minimal\/} and $s$-{\it maximal\/} ones),
establish bijections between these two sets of elements and the coroot
lattice points of certain simplices, and  eventually obtain formulae for
the number of short and strictly $s$-positive antichains.
Second, we introduce
and study the {\it semi-Catalan arrangement\/}, which has the same relation
to short and strictly $s$-positive antichains as the usual Catalan
arrangement has to all and strictly positive antichains.
The difference between the Catalan and semi-Catalan arrangements is that
we "deform" only the hyperplanes orthogonal to short roots in the latter.
We prove various results connecting the
dominant regions of the semi-Catalan
arrangement and the elements of $\HW$ attached to short antichains.
Adapting Athanasiadis' argument from \cite{ath02}, we compute the characteristic
polynomial  for the {\it extended semi-Catalan arrangements\/},
or in other words, for
$m$-semi-Catalan arrangements, $\catsm$, with $m=0,1,2,\dots$.
For $m=0$, one obtains the Coxeter arrangement of $W$, and for $m=1$,
the semi-Catalan arrangement.
\\[.4ex] \indent
Here is a part of our results. Let $\ap_1,\dots,\ap_p$ be the simple roots
of $\g$ and  $\theta$ the highest root.
Let $\gA$ be the fundamental alcove of $\HW$ and
$g$ the sum of coefficients of the {\sl short} simple roots in the
expression of $\theta=\sum c_i\ap_i$.
Then the short (resp. strictly
$s$-positive) antichains are in a one-to-one correspondence with
the coroot lattice points in $(g+1)\gA$ (resp. $(g-1)\gA$).
If the root system is not of type $\GR{G}{2}$, this leads to a closed
formula for the number of the respective antichains. E.g.,
the number of short antichains
is equal to $\displaystyle\prod_{i=1}^p \frac{g+e_i+1}{e_i+1}$,
where $e_i$, $i=1,2,\dots,p$, are the exponents of the Weyl group $W$.
Using this, we found a uniform expression, which covers the
$\GR{G}{2}$-case as well, see Eq.~\re{unified},
but it awaits for a conceptual explanation. The characteristic
polynomial of $\catsm$ is
$\chi(t)=\prod_{i=1}^p(t-mg-e_i)$ (again, if $\Delta$ is not of type
$\GR{G}{2}$.) For $\GR{G}{2}$, the formula for
$\chi(t)$ depends on the parity of $m$.
We also define a `short' analogue of the extended Shi arrangement, which we
call, of course, the {\it extended semi-Shi arrangement}, and propose
a conjectural formula for its characteristic polynomial, see \re{semi-shi}.
\\[.4ex] \indent
A rough description of the contents is as follows.
In Sections~\ref{old} and \ref{reg}, we give a review of results
concerning ideals (antichains) and Catalan
arrangements, including two approaches described above.
In particular, we consider minimal and maximal elements of
$\HW$ and their connection with ideals.
Some complements to known results are also given.
We attempt
to present a unified treatment that can be generalized afterwards,
without much pains, to the setting of short antichains. Our main
results are gathered in Sections~\ref{short}-\ref{numer}.
After a brief description in Section~\ref{short} of the relationship between
$\be$-stable subspaces of the little adjoint $G$-module and short
antichains, we turn, in Section~\ref{new}, to considering $s$-minimal and
$s$-maximal elements of $\HW$ and related simplices.
In Section~\ref{geners}, we compute the characteristic polynomial
for the $m$-semi-Catalan arrangement with arbitrary $m\in\Bbb N$ and
study relationship between the semi-Catalan arrangement (which corresponds
to $m=1$) and short antichains.
As a consequence of our theory, we present, in Section~\ref{numer},
several intriguing results whose proof uses case-by-case verification.
\\[.4ex] \indent
To a great extent, this  work was inspired by recent preprints
of \ C.\,Athanasiadis \cite{ath02} and E.\,Sommers \cite{eric}.

{\small {\bf Acknowledgements.} This paper was written during my stay at the
Max-Planck-Institut f\"ur Mathematik (Bonn).
I would like to thank this institution for hospitality
and excellent working conditions.
}

\section{Notation and other preliminaries}
\label{prelim}

\noindent
\begin{subs}{Main notation}
\end{subs}
$\Delta$ is the root system of $(\g,\te)$ and
$W$ is the usual Weyl group. For $\ap\in\Delta$, $\g_\ap$ is the
corresponding root space in $\g$.

$\Delta^+$  is the set of positive
roots and $\rho=\frac{1}{2}\sum_{\ap\in\Delta^+}\ap$.

$\Pi=\{\ap_1,\dots,\ap_p\}$ is the set of simple roots in $\Delta^+$
and $\theta$ is the highest root in $\Delta^+$.
 \\
 We set $V:=\te_{\Bbb R}=\oplus_{i=1}^p{\Bbb R}\ap_i$ and denote by
$(\ ,\ )$ a $W$-invariant inner product on $V$. As usual,
$\mu^\vee=2\mu/(\mu,\mu)$ is the coroot
for $\mu\in \Delta$.

${\gC}=\{x\in V\mid (x,\ap)>0 \ \ \forall \ap\in\Pi\}$
\ is the (open) fundamental Weyl chamber.

${\gA}=\{x\in V\mid (x,\ap)>0 \ \  \forall \ap\in\Pi \ \ \& \
(x,\theta)<1\}$ \ is the fundamental alcove.


$Q^+=\{\sum_{i=1}^p n_i\ap_i \mid n_i=0,1,2,\dots \}$
and $Q^\vee=\oplus _{i=1}^p {\Bbb Z}\ap_i^\vee\subset V$
is the coroot lattice.
\\
Letting $\widehat V=V\oplus {\Bbb R}\delta\oplus {\Bbb R}\lb$, we extend
the inner product $(\ ,\ )$ on $\widehat V$ so that $(\delta,V)=(\lb,V)=
(\delta,\delta)=
(\lb,\lb)=0$ and $(\delta,\lb)=1$.

$\widehat\Delta=\{\Delta+k\delta \mid k\in {\Bbb Z}\}$ is the set of affine
real roots and $\widehat W$ is the  affine Weyl group.
\\
Then $\widehat\Delta^+= \Delta^+ \cup \{ \Delta +k\delta \mid k\ge 1\}$ is
the set of positive
affine roots and $\widehat \Pi=\Pi\cup\{\ap_0\}$ is the corresponding set
of affine simple roots,
where $\ap_0=\delta-\theta$.
The inner product $(\ ,\ )$ on $\widehat V$ is
$\widehat W$-invariant. The notation $\beta>0$ (resp. $\beta <0$)
is a shorthand for $\beta\in\HD^+$ (resp. $\beta\in -\HD^+$).
\\
For $\ap_i$ ($0\le i\le p$), we let $s_i$ denote the corresponding simple
reflection in $\widehat W$.
If the index of $\ap\in\widehat\Pi$ is not specified, then we merely write
$s_\ap$. 
The length function on $\widehat W$ with respect
to  $s_0,s_1,\dots,s_p$ is denoted by $\ell$.
For any $w\in\widehat W$, we set
\[
   N(w)=\{\ap\in\widehat\Delta^+ \mid w(\ap) \in -\widehat \Delta^+ \} .
\]
It is standard that $\#N(w)=\ell(w)$ and $N(w)$ is {\it bi-convex\/}.
The latter means
that both $N(w)$ and $\HD^+\setminus N(w)$ are subsets of $\HD^+$
that are closed under addition.
Furthermore, the assingment $w\mapsto N(w)$ sets up
a bijection between the elements of $\HW$ and the finite bi-convex subsets
of $\HD^+$.

\begin{subs}{Ideals and antichains}               \label{ideals}
\end{subs}
Throughout the paper, $\be$ is the Borel subalgebra of $\g$ corresponding
to $\Delta^+$ and $\ut=[\be,\be]$.
Let $\ce\subset\be$ be an {\sf ad}-nilpotent ideal. Then
$\ce=\underset{\ap\in I}{\oplus}\g_\ap$
for some $I\subset \Delta^+$. This $I$ is said to be an {\it ideal\/} (of
$\Delta^+$). More precisely, a set $I\subset\Delta^+$ is an ideal, if
whenever $\gamma\in I,\mu\in\Delta^+$, and $\gamma+\mu\in\Delta$, then
$\gamma+\mu\in I$.
Our exposition will be mostly combinatorial, i.e., in place of
{\sf ad}-nilpotent ideal of $\be$ we will deal with the respective
ideals of $\Delta^+$.
\\[.5ex]
For $\mu,\gamma\in\Delta^+$, write $\mu\curle\gamma$, if
$\gamma-\mu\in Q^+$. The notation $\mu\prec\gamma$ means that
$\mu\curle\gamma$ and $\gamma\ne\mu$.
We regard $\Delta^+$ as poset under `$\curle$'.
Let $I\subset\Delta^+$ be an ideal. An element $\gamma\in I$ is called
a {\it generator\/}, if $\gamma-\ap\not\in I$ for any $\ap\in\Pi$.
In other words, $\gamma$ is a minimal element of $I$ with respect to
`$\curle$'.
We write $\Gamma(I)$ for the set of generators of $I$.
It is easily seen that $\Gamma(I)$ is an {\it antichain\/}
of $\Delta^+$, i.e., $\gamma_i\not\curle\gamma_j$ for any pair
$(\gamma_i,\gamma_j)$ in $\Gamma(I)$.
Conversely, if $\Gamma\subset \Delta^+$ is an antichain,
then the ideal
\[
I\langle\Gamma\rangle:=\{\mu\in \Delta^+\mid \mu\succcurlyeq\gamma_i
\text{ for some } \gamma_i\in \Gamma \}
\]
has $\Gamma$ as the set of generators.
Let $\AN$ denote the set of all antichains in $\Delta^+$.
In view of the above bijection
$\AD \overset{1:1}{\longleftrightarrow}\AN$,
we will freely switch between
ideals and antichains. An ideal $I$ is called {\it strictly positive\/}, if
$I\cap\Pi=\varnothing$. The set of strictly positive ideals is denoted
by $\AD_0$.


\section{Ideals, maximal and minimal elements of $\HW$}
\label{old}
\setcounter{equation}{0}

\noindent
In this section we  review some recent results by Athanasiadis,
Cellini-Papi, Sommers, and this author. A few complements is also given.
\\[.5ex]
The idea of describing ideals of $\Delta^+$ through the use of elements of
$\HW$ goes back to D.~Peterson, who exploited minuscule elements for
counting Abelian ideals of $\be$, see \cite{Ko1}.
In general case, given $I\subset \Delta^+$, we want to have $w\in\HW$
such that $N(w)\subset \cup_{k\ge 1} (k\delta-\Delta^+)$ and
$N(w)\cap (\delta-\Delta^+)=\delta-I$. It turns out that, for any ideal
$I$, there is a unique element of minimal length satisfying these
properties. In contrast,
the element of maximal length exists if and only if $I$ is strictly
positive, and in this case such an element is unique, too.
Implementation of this program yields also explicit formulae for
the number of all and strictly positive ideals.

As is well known, $\widehat W$ is isomorphic to a semi-direct product
of $W$ and $Q^\vee$.
Given $w\in\widehat W$, there is a unique decomposition
\begin{equation}  \label{affine}
w=v{\cdot}t_{r}\ ,
\end{equation}
where $v\in W$ and  $t_{r}$ is the translation
corresponding to $r\in Q^\vee$.
The word ``translation" means the following.  The group $\HW$ has two
natural actions:

(a) the linear action on $\widehat V=
V\oplus{\Bbb R}\delta\oplus{\Bbb R}\lb$;

(b) the affine-linear action on $V$.
\\
We use `$\ast$' to denote the second action.
For $r\in Q^\vee$, the linear action of $t_r\in \widehat W$ on
$V\oplus {\Bbb R}\delta$ is given by $t_r(x)=x-(x,r)\delta$ (we do not need
the formulas
for the whole of $\widehat V$), while the affine-linear action on $V$ is
given by $t_r\asst y=y+r$. So that $t_r$ is a true translation for the
$\ast$-action on $V$.
\\[.5ex]
Let us say that $w\in\HW$ is {\it dominant\/}, if
$w(\ap)>0$ for all $\ap\in\Pi$.
Obviously, $w$ is dominant if and only if
$N(w)\subset \cup_{k\ge 1} (k\delta-\Delta^+)$.
It also follows from \cite[1.1]{CP1} that
$w$ is dominant if and only if
$w^{-1}\asst \gA\subset \gC$. Write $\HW_{dom}$ for the set of
dominant elements.

\begin{s}{Proposition}  \label{opis-dom}
\begin{itemize}
\item[\sf (i)] If $w=v{\cdot}t_r\in \HW_{dom}$, then $r\in -\ov{\gC}$;
\item[\sf (ii)] The mapping $\HW_{dom}\to Q^\vee$ given
by $w=v{\cdot}t_r\mapsto v(r)$ is a bijection.
\end{itemize}
\end{s}\begin{proof}
(i) \ We have $w^{-1}\asst x=v^{-1}(x)-r$ for any $x\in V$.
In particular, $w^{-1}\asst 0=-r$.
Since $0\in\ov{\gA}$ and $w$ is dominant, we are done.

(ii) \ Given $\varkappa\in Q^\vee$, we want
to find $w=v{\cdot}t_r$ such that
$w^{-1}\asst\gA=v^{-1}(\gA)-r\subset \gC$ and $v(r)=\varkappa$.
In view of the last equality, the previous containment reads
$v^{-1}(\gA-\varkappa)\subset \gC$. Therefore $v$ must be the
unique element
of $W$ taking the alcove $\gA-\varkappa$ into the dominant Weyl chamber
$\gC$. Then $r=v^{-1}(\varkappa)$.
\\
This argument proves both the injectivity and surjectivity of the
mapping in question.
\end{proof}%
Letting $\delta-I:=N(w)\cap (\delta-\Delta^+)$, we easily
deduce that $I$ is an ideal, if $w\in \HW_{dom}$. We say $\delta-I$ is
the {\it first layer\/} of $N(w)$ and $I$ is the
{\it first layer ideal\/} of $w$.
However, an ideal $I$ may well arise from different dominant elements.
To obtain a bijection, one has to impose further constraints on
dominant elements. One may attempt to consider either maximal or minimal
bi-convex subsets with first layer $\delta-I$.
This naturally leads to notions of
`minimal' and `maximal' elements. This terminology suggested in
\cite{eric} is also explained by a relationship between these elements
and dominant regions of the Shi arrangement, see Section~\ref{reg}.
However, the formal definitions
do not require invoking arrangements. Furthermore, we want to stress
that many results relating the ideals and these two kinds of dominant
elements can be obtained without ever mentioning the Shi (or Catalan)
arrangement.

\begin{rem}{Definition}   \label{def-mi}
$w\in\HW$ is called {\it minimal\/}, if
\begin{itemize}
\item[\sf (i)] \ $w$ is dominant;
\item[\sf (ii)] \ if $\ap\in\HP$ and
$w^{-1}(\ap)=k\delta+\mu$ for some $\mu\in \Delta$, then $k\ge -1$.
\end{itemize}
\end{rem}%
Using (i), condition (ii) can be made more precise. If $k\in\{-1,0\}$,
then $\mu\in \Delta^+$.
\\
The set of minimal elements is denoted by $\HW_{min}$.

\begin{s}{Proposition {\ququ \cite[Prop.\,2.12]{CP1}}}
\label{bij-mi}
There is a bijection between $\HW_{min}$ and\/ $\AD$.
Namely,
\begin{itemize}
\item given $w\in \HW_{min}$, the corresponding ideal is
$\{\mu\in\Delta^+\mid \delta-\mu\in N(w)\}$;
\item given $I\in\AD$, the corresponding minimal element is determined by
the finite bi-convex set
\[
   \bigcup_{k\ge 1}(k\delta-I^k)\subset \HD^+ \ .
\]
Here $I^k$ is defined inductively by $I^k=(I^{k-1}+I)\cap\Delta^+$.
\end{itemize}
\end{s}%
If $N\subset \HD^+$ is a finite convex subset, containing $\delta-I$,
then it must also contain $\cup_{i\ge 1}(k\delta-I^k)$. So, the
latter is the minimal bi-convex subset containing $\delta-I$.
\\[.5ex]
The first layer ideal of $w\in\HW_{min}$ is denoted by $I_w$.

\begin{s}{Proposition {\ququ \cite[Theorem\,2.2]{duality}
\cite[6.3(1)]{eric}}}
\label{gen-mi}  \par
If $w\in\HW_{min}$, then
$\Gamma(I_w)=\{\gamma\in\Delta^+ \mid w(\delta-\gamma) \in-\HP\}$.
\end{s}%
Following \cite{CP2}, we give a "geometric" description of the
minimal elements.
Set
\[
 D_{min}=\{x\in V \mid (x,\ap)\ge -1 \ \forall\ap\in\Pi \ \ \& \
\ (x,\theta)\le 2\} \ .
\]
It is a certain simplex in $V$.

\begin{s}{Proposition {\ququ \cite[Prop.\,2 \& 3]{CP2}}}
\label{opis-mi} \par
1. $w=v{\cdot}t_r\in \HW_{min}$ \ $\Longleftrightarrow
\left\{ \begin{array}{l} w \text{ is dominant},  \\
                         v(r)\in D_{min}\cap {Q}^\vee \ .
       \end{array}\right.$

2. The mapping $\HW_{min} \to D_{min}\cap Q^\vee$,
$w=v{\cdot}t_r\mapsto v(r)$, is a bijection.
\end{s}\begin{proof}
1. `$\Rightarrow$'
The first condition is satisfied by the definition. \\
Next, we have $w^{-1}(x)=v^{-1}(x)+(x, v(r))\delta$ \ for any
$x\in V\oplus {\Bbb R}\delta$. In particular,
\begin{equation} \label{simroots} \begin{array}{l}
  w^{-1}(\ap_i)=v^{-1}(\ap_i)+(\ap_i, v(r))\delta, \quad i\ge 1 ,\\
  w^{-1}(\ap_0)=-v^{-1}(\theta)+(1-(\theta, v(r)))\delta \ .
\end{array}
\end{equation}
Comparing this with Definition~\ref{def-mi}(ii),
one concludes that $v(r)\in D_{min}$.

`$\Leftarrow$'  The previous argument can be reversed.

2. This follows from part~1 and Proposition~\ref{opis-dom}.
\end{proof}%
{\sf Remark}. The above proof applies equally well to
Propositions~\ref{opis-ma},\,\ref{opis-smi} and
\ref{opis-sma} below. It is a simplified version
of the proof of Propositions~2 \& 3 in \cite{CP2}.
\\[.5ex]
It follows that $\#(\AD)$ equals the number
of integral points in $D_{min}$. (Unless otherwise stated, an
'integral point' is a point lying in $Q^\vee$.)
A pleasant feature of this situation is that there
is an element of $\HW$ that takes $D_{min}$ to a dilated
closed fundamental alcove. Namely,
$w(D_{min})=(h+1)\ov{\gA}$ for some $w\in \HW$, see \cite[Thm.\,1]{CP2}.
Write $\theta$ as a linear combination of simple roots:
$\theta= \sum_i c_i\ap_i$. The integers $c_i$ are said to be the
{\it coordinates of\/} $\theta$.
By a result of M.\,Haiman \cite[7.4]{mark}, the number of integral points
in $t\ov{\gA}$ is equal to
\begin{equation}  \label{haiman}
\prod_{i=1}^p \frac{t+e_i}{1+e_i} \
\end{equation}
whenever $t$ is relatively prime with all the coordinates
of $\theta$. Since this condition is satisfied for
$t=h+1$, one obtains
\begin{equation}  \label{chislo-mi}
   \# (\AD)=\prod_{i=1}^p \frac{h+e_i+1}{e_i+1} \ .
\end{equation}
It is the main result of \cite{CP2}.
\\
Combining Propositions~\ref{gen-mi} and Eq.~\re{simroots} yields the
assertion that $\#\Gamma(I_w)=k$ if and only if $v(r)$ lies on a face of
$D_{min}$ of codimension $k$ \cite[Thm.\,2.9]{duality}.
\vskip.6ex\noindent
Now, we turn to considering  maximal (dominant) elements of $\HW$.
Most of the results on these elements are due to E.~Sommers.
However our presentation follows \cite{eric} rather freely,
since we want to have a uniform treatment for both
minimal and maximal elements. Because some
assertions have no exact counterparts in \cite{eric},
we also give some proofs.

\begin{rem}{Definition}   \label{def-ma}
$w\in\HW$ is called {\it maximal\/}, if
\begin{itemize}
\item[\sf (i)] $w$ is dominant;
\item[\sf (ii)] if $\ap\in\HP$ and
$w^{-1}(\ap)=k\delta+\mu$ for some $\mu\in \Delta$, then $k\le 1$.
\end{itemize}
\end{rem}%
Using (i), condition (ii) can be made more precise.
If $k=1$, then $\mu\in -\Delta^+$; if $k=0$, then $\mu\in \Delta^+$
The set of maximal elements is denoted by
$\HW_{max}$.

If $I\in\AD_0$, then for any $\mu\in\Delta^+$ we define
$k(\mu, I)$ as the minimal possible number of summands
in the expression $\mu=\sum_{i}\nu_i$, where $\nu_i\in \Delta^+\setminus I$.
Notice that this definition only makes sense for strictly positive ideals.

\begin{s}{Proposition {\ququ \cite[Sect.\,5]{eric}}}
\label{bij-ma}
There is a bijection between $\HW_{max}$ and $\AD_0$.
Namely,
\begin{itemize}
\item \ given $w\in \HW_{min}$, the corresponding strictly
positive ideal is
$\{\mu\in\Delta^+\mid \delta-\mu\in N(w)\}$;
\item \ given $I\in\AD_0$, the corresponding maximal element is
determined by the finite bi-convex set
\end{itemize}
\hbox to \textwidth{\enspace ($\Diamond$)\hfil
 $\{ m\delta-\gamma\mid \gamma\in I \quad\&\quad
 1\le m\le k(\gamma,I){-}1\}$. \hfil}
\end{s}\begin{proof}
1. Suppose $w\in\HW$ is dominant, and
let $I$ be the first layer ideal of $w$.
Assuming that $I\cap\Pi\ni \ap$, we show that $w$
cannot be maximal.
For any $\gamma\in I$, let $k_\gamma$ be the maximal integer such that
$k_\gamma\delta-\gamma\in N(w)$, i.e.,
\[
  N(w)=\{l\delta-\gamma\mid \gamma\in I \ \& \ 1\le l\le k_\gamma\} \ .
\]
Let $I(\ap)$ be the ideal generated by $\ap$. Clearly, $I(\ap)\subset I$.
Set
\[
  N(w)^{\langle 2\rangle}=
\{l\delta-\gamma\mid \gamma\in I(\ap) \ \& \ 1\le l\le 2k_\gamma\}\cup
\{l\delta-\gamma\mid \gamma\in I\setminus I(\ap)
\ \& \ 1\le l\le k_\gamma\} \ .
\]
Obviously, $N(w)^{\langle 2\rangle}$ is finite and has the same first layer
as $N(w)$. It is also easy to verify that $N^{\langle 2\rangle}$ is again
bi-convex.
Hence  $N(w)^{\langle 2\rangle}=N(w')$
for some $w'\in\HW$. Since $N(w')\supset N(w)$, there is a
presentation $w'=uw$, where $\ell(w')=\ell(u)+\ell(w)$. If $s_\nu$ ($\nu\in\HP$)
is the rightmost
reflection in a reduced decomposition for $u$, then
$w^{-1}(\nu)=k\delta-\mu$ with $k\ge 2$, as the first layers of
$N(w')$ and $N(w)$ are the same. Thus, $w$ is not maximal.

2. Suppose $I\in\AD_0$, and let $w\in\HW$ be any dominant element with
first layer ideal $I$. Since $\HD^+\setminus N(w)$ is convex and
contains $\delta-(\Delta^+\setminus I)$, it follows
from the very definition of numbers
$k(\gamma,I)$ that $l\delta-\gamma\in \HD^+\setminus N(w)$
for all $l\ge k(\gamma,I)$. Hence $N(w)$ is contained in the finite
set given by Eq.~\ref{bij-ma}($\Diamond$). It only remains to prove
that the latter is bi-convex. For this crucial fact, we refer to
\cite[Lemma\,5.2]{eric}.
\end{proof}%
The strictly positive ideal corresponding to $w\in\HW_{max}$
(the first layer ideal of $w$) is denoted
by $I^w$.
For an ideal $I\subset\Delta^+$, we write $\Xi(I)$ for the set of
maximal elements of $\Delta^+\setminus I$. It is immediate that $\Xi(I)$
is an antichain.

\begin{s}{Proposition {\ququ \cite[6.3(2)]{eric}}}  \label{gen-ma}
If $w\in\HW_{max}$, then
$\Xi(I^w)=\{\gamma\in\Delta^+\mid w(\delta-\gamma)\in\HP\}$.
\end{s}%
\vskip-1.2ex
\begin{rem}{Remark}    \label{cover}
Note that antichains of the form $\Xi(I^w)$ are not arbitrary. From
the definition of a strictly positive ideal it readily follows that,
given $\Xi\in \AN$, we have
$\Xi=\Xi(I)$ for some $I\in \AD_0$ if and only if for any
$\ap\in\Pi$ there is a $\gamma\in\Xi$ such that $\gamma\succcurlyeq \ap$.
We shall say that such an antichain {\it covers\/} the simple roots.
\end{rem}%
Now, we proceed to a "geometric" characterization of the maximal elements.
Set
\[
 D_{max}=\{x\in V \mid (x,\ap)\le 1 \ \forall\ap\in\Pi \ \ \& \ \
(x,\theta)\ge 0\} \ .
\]
It is a certain simplex in $V$.

\begin{s}{Proposition {\ququ (cf. \cite[Prop.\,5.6]{eric})}}
\label{opis-ma} \par
1. $w=v{\cdot}t_r\in \HW_{max}$ \ $\Longleftrightarrow
\left\{ \begin{array}{l} w \text{ is dominant},  \\
                         v(r)\in D_{max}\cap {Q}^\vee \ .
       \end{array}\right.$

2. The mapping $\HW_{max} \to D_{max}\cap Q^\vee$,
$w=v{\cdot}t_r\mapsto v(r)$, is a bijection.
\end{s}\begin{proof*}
1. The argument is the same as in Proposition~\ref{opis-mi}, taking into
account that the constraints for $D_{max}$ are different.
\\
2. This follows from part~1 and Proposition~\ref{opis-dom}.
\end{proof*}%

\begin{s}{Proposition}  \label{cor-opis-ma}
Suppose $w=v{\cdot}t_r\in \HW_{max}$. Then $\#\Xi(I^w)=k$ if and only
$v(r)$ lies on a face of codimension $k$ of $D_{max}$.
\end{s}\begin{proof}
Combine Propositions~\ref{gen-ma} and Eq.~\re{simroots}.
\end{proof}%
Since $\Pi$ is the only antichain of cardinality $p$
\cite[2.10(ii)]{duality} and it is
certainly of the form $\Xi(I^w)$, we see that $D_{max}$ has a unique
integral vertex.
\\[.5ex]
In order to compute  $\#(D_{max}\cap Q^\vee)$,
we replace $D_{max}$ with another simplex.
Let $\{\varpi^\vee_i\}_{i=1}^p$ denote the dual basis of $V$ for
$\{\ap_i\}_{i=1}^p$. Set $\rho^\vee=\sum_{i=1}^p\varpi^\vee_i$.
Since the sum of the coordinates of $\theta$ equals $h-1$, the translation
$x\mapsto t_{-\rho^\vee}\asst x=x-\sum_{i=1}^p\varpi^\vee_i$
takes $D_{max}$ to the negative dilated fundamental alcove
\[
-(h-1)\ov{\gA}=\{ x\in V\mid (x,\ap)\le 0 \quad \forall\,\ap\in\Pi;\
   (x,\theta)\ge 1-h\} .
\]
It may happen that $\rho^\vee$ does not belong to $Q^\vee$, so that
this translation, which is in
the extended affine Weyl group, does not belong to $\HW$,
while we wish to have a transformation from $\HW$.
Nevertheless, since $h-1$ is relatively prime with the index of
connection of $\Delta$, it follows from \cite[Lemma\,1]{CP2} that
there is an element of
$\HW$ that takes $D_{max}$ to $(1-h)\ov{\gA}$.

Again, using the above-mentioned result of Haiman, see Eq.~\re{haiman},
one obtains

\begin{s}{Theorem {\ququ \cite{ath02},\cite{duality},\cite{eric}}}
\label{chislo-ma}
\[
   \#(\AD_0)= \prod_{i=1}^p \frac{h+e_i-1}{e_i+1}  \ .
\]   \vskip-1ex
\end{s}%
{\sf Remark.} The proofs in \cite{ath02} and \cite{duality} are based on
the fact that the strictly positive ideals correspond to the bounded
regions of the Catalan arrangement and that the number of bounded regions
of any hyperplane arrangement can be computed  via the characteristic
polynomial of this arrangement, see Section~\ref{reg}.


\section{Ideals and dominant regions of the Catalan arrangement}
\label{reg}
\setcounter{equation}{0}

\noindent
Recall a bijection between the ideals of $\Delta^+$
and the dominant regions of the Catalan arrangement.
This bijection is due to
Shi \cite[Theorem\,1.4]{shi}.  
\\[.5ex]
For $\mu\in\Delta^+$ and $k\in{\Bbb Z}$, define the hyperplane
$\gH_{\mu,k}$ in $V$ as $\{x\in V \mid (x,\mu)=k\}$.
The {\it Catalan arrangement\/}, $\mathrm{Cat}(\Delta)$,
is the collection of hyperplanes
$\gH_{\mu,k}$, where $\mu\in \Delta^+$ and $k=-1,0,1$.
The {\it regions\/} of an arrangement are the connected components of
the complement in $V$ of the union of all its hyperplanes.
Obviously, the dominant regions of $\mathrm{Cat}(\Delta)$ are the same as
those for the {\it Shi arrangement\/} $\mathrm{Shi}(\Delta)$.
The latter is the collection of hyperplanes
$\gH_{\mu,k}$, where $\mu\in \Delta^+$ and $k=0,1$.
But, it will be more convenient for us to deal with the arrangement
$\mathrm{Cat}(\Delta)$, since it is $W$-invariant.
\\[.5ex]
It is clear that $\gC$ is a union of regions of $\mathrm{Cat}(\Delta)$.
Any region lying in $\gC$ is said to be
{\it dominant\/}.
The Shi bijection takes an ideal $I\subset \Delta^+$ to the dominant region
\begin{equation}   \label{bij-shi}
 R_I =\{ x\in \gC \mid  (x,\gamma)>1, \text{ if \ }\gamma\in I \quad \&\quad
           (x,\gamma)<1,  \text{ if \ }\gamma\not\in I
           \}\ .
\end{equation}
A proof of this result, which is almost self-contained,
is given by Athanasiadis in \cite[Lemma\,3.1]{ath03}.
But for the fact that $R_I$ is non-empty he refers to an earlier
work Shi \cite{shi0}. It is not, however, easy to extract the proof
from Shi's paper. What we want to say is that the non-emptiness of
$R_I$ readily follows from the theory of minimal elements:

{\it If $w\in\HW$ is the minimal element corresponding to $I$, then
$w^{-1}\asst\gA$ belong to $R_I$.}

\noindent
Indeed, $\gH_{\mu,1}$ separates $\gA$ and $w^{-1}\asst\gA$ if and only
if $w(\delta-\mu)\in -\HD^+$, see \cite[1.1]{CP1}.
\\
In fact, $w^{-1}\asst\gA$ is the nearest to the origin alcove in $R_I$.

A region (of an arrangement) is called {\it bounded\/}, if it is contained
in a sphere about the origin.

\begin{s}{Proposition {\ququ \cite{ath02},\cite{duality},\cite{eric}}}
\label{regions}
$I\in\AD(\g)_0$
if and only if the region $R_I$ is bounded.
\end{s}%
If $R_I$ is bounded, then it obviously contains an alcove that is
most distant from the origin. It was shown in \cite{eric} that if
$w$ is the maximal element corresponding to $I\in\AD_0$, then
$w^{-1}\asst\gA$ is the most distant from the origin alcove in $R_I$.
\\[.5ex]
The number of regions and bounded regions of any hyperplane arrangement
can be counted through the use of a striking result of T.\,Zaslavsky.
Let $\chi(\ca, t)$ denote the characteristic polynomial of a
hyperplane arrangement $\ca$ in $V$ (see e.g. \cite[Sect.\,2]{ath02} for
precise definitions).
\begin{s}{Theorem {\ququ (Zaslavsky)}}
\label{zaslavsky}   \par
1. The number of regions into which $\ca$ dissects $V$ equals
$(-1)^p\chi(\ca,-1)$.
\par 
2. The number of bounded regions into which $\ca$ dissects $V$ equals
$|\chi(\ca,1)|$.
\end{s}%
In \cite{ath02}, Athanasiadis gives a nice case-free proof of
the following formula for the characteristic polynomial of
the Catalan arrangement:
\begin{equation} \label{char-cat}
  \chi(\mathrm{Cat}(\Delta),t)=\prod_{i=1}^p (t-h-e_i) \ .
\end{equation}
Since $\mathrm{Cat}(\Delta)$ is $W$-invariant,
the values $\displaystyle
\frac{|\chi(\mathrm{Cat}(\Delta),\pm 1)|}{\#(W)}$ give
the number of bounded and all regions in $\gC$, respectively.
In this way, one obtains explicit formulae for the cardinality
of $\AD_0$ and $\AD$ written already down in Section~\ref{old}.
Thus, the characteristic polynomial of the
Catalan arrangement provides an alternative approach to counting
ideals and strictly positive ideals.


\section{Short antichains and $\be$-stable subspaces in the little
adjoint $G$-module}
\label{short}
\setcounter{equation}{0}

\noindent
At the rest of the paper, we stick to the case in which $\Delta$ has
roots of different length. Then we naturally have long and short
roots, long and short reflections, etc. Our goal is to show that
the theory presented in the previous sections can be extended to
the setting, where one pays attention to the length of roots involved.
A piece of such theory already appeared in \cite{ablong}, where
we studied Abelian ideals of $\Delta^+$ consisting
of only long roots. Now, we consider the general case.
Our treatment will again be combinatorial. We wish, however, to stress
that it has a related representation-theoretic picture.
While the ideals (antichains) in $\Delta^+$ correspond bijectively
to the $\be$-stable subspaces in $\g$ having no nonzero semisimple
elements, our short antichains in $\Delta^+$ correspond
bijectively to the $\be$-stable subspaces, without
nonzero semisimple elements, in the little adjoint $\g$-module.
\\[.6ex]
To distinguish various objects associated with
long and short roots, we use the subscripts `$l$' and `$s$', respectively.
For instance, $\Pi_l$ is the set of long simple roots
and $\Delta^+_s$ is the set of short positive roots. Accordingly, each
simple reflection $s_i$ is either short or long. Since $\theta$ is long,
the simple root $\ap_0$ and the reflection $s_0$ are regarded as long.
Therefore, $\HP_l=\Pi_l\cup\{\ap_0\}$.
Write $\theta_s$ for the unique short dominant root in $\Delta^+$.
A simple $\g$-module with highest weight $\theta_s$, $\vts$, is said to be
{\it little adjoint}.
The set of nonzero weights of $\vts$ is $\Delta_s$, all nonzero weights
are simple, and the multiplicity of the zero weight is $\#(\Pi_s)$
\cite[2.8]{ya-tg}.

\begin{rem}{Definition}  \label{def-shan}
An antichain $\Gamma\subset\Delta^+$ is called {\it short}, if it
consists of short roots, i.e., $\Gamma\subset\Delta^+_s$. Similarly,
one defines a {\it long\/} antichain.
\end{rem}%
If $\Gamma$ is a short antichain, then $\Gamma^\vee$ is a long
antichain in the dual root system $\Delta^\vee$. Therefore, it suffices,
in principle, to consider only short antichains. We write $\AN_s$
for the set of all short antichains of $\Delta^+$.  The respective set
of ideals is denoted by $\AD_s$.
\\
Recall that, for any finite-dimensional rational $G$-module $\VV$, there
are notions of semisimple
and nilpotent elements, generalizing those in $\g$, see \cite{VP}.
An element  $v\in\VV$ is called {\it semisimple},  if
the orbit $Gv$ is closed; it is called {\it nilpotent}, if
the closure of $Gv$ contains the origin. We shall say that a subspace
of $\VV$ is {\it nilpotent\/}, if it consists of nilpotent elements.

\begin{s}{Proposition}   \label{no-ss}
There is a one-to-one correspondence between $\AN_s$ and
the nilpotent $\be$-stable subspaces of
$\vts$.
\end{s}\begin{proof}
If $\Gamma$  is a short antichain, then the corresponding
subspace of $\vts$ is the sum of weight spaces $\vts_\mu$,
$\mu\in\Delta^+_s$,
such that $\gamma_i\cyeq \mu$ for some $\gamma_i\in\Gamma$.
Conversely, if $\UU$ is a $\be$-stable subspace of $\vts$, then
it is a sum of weight spaces. If $\UU$ contains no semisimple elements,
then all weights of $\UU$ must lie in an open halfspace of $V$
(by the Hilbert-Mumford criterion). From
the above description of weights of $\vts$, it follows that
the weights of $\UU$ must form a subset of $\Delta^+_s$.
The minimal elements of this set of weights give us the required
short antichain.
\end{proof}%
If $\Gamma$ is a short antichain, then the set of weights of the
corresponding nilpotent
$\be$-stable subspace of $\vts$ is $I\langle\Gamma\rangle\cap\Delta^+_s$.


\section{Short antichains, $s$-minimal and $s$-maximal elements of $\HW$}
\label{new}
\setcounter{equation}{0}

\noindent
Our goal is to show that the theory described in Section~\ref{old}
extends well to short antichains.


\begin{rem}{Definition}   \label{def-smi}
$w\in\HW$ is called {\it $s$-minimal\/}, if
\begin{itemize}
\item[\sf (i)] $w$ is dominant;
\item[\sf (ii)]
if $\ap\in\Pi_s$ and $w^{-1}(\ap)=k\delta+\mu$ with
$\mu\in \Delta$, then $k\ge -1$;
\item[\sf (iii)]
if $\ap\in\HP_l$ and $w^{-1}(\ap)=k\delta+\mu$ with
$\mu\in \Delta$, then $k\ge 0$;
\end{itemize}
\end{rem}%
Using (i), conditions (ii),\,(iii) can be made more precise.
If $k=0$ or $k=-1$ in (ii), then $\mu\in \Delta^+$.

\noindent
We write $\HW^{(s)}_{min}$ for the set of all $s$-minimal
elements. Notice that $\HW^{(s)}_{min}\subset \HW_{min}$.

\begin{s}{Proposition}
The bijection between $\AD$ and $\HW_{min}$ described in
Proposition~\ref{bij-mi} gives rise to a bijection between
$\AN_s$ (or $\AD_s$) and $\HW^{(s)}_{min}$.
\end{s}\begin{proof}
1. Suppose $w\in\HW^{(s)}_{min}$, and let $I_w$ be the corresponding
ideal. It follows from Definition~\ref{def-smi}(ii),(iii) and
Proposition~\ref{gen-mi} that $\Gamma(I_w)\subset \Delta^+_s$.
Thus, we obtain a short antichain.

2. The use of Proposition~\ref{gen-mi} gives also the converse.
\end{proof}%
Now, we give a geometric description of $s$-minimal elements
in the spirit of Section~\ref{old}.
Set
\[
  D^{(s)}_{min}=\{x\in V \mid (x,\ap)\ge -1 \ (\ap\in\Pi_s); \quad
(x,\ap)\ge 0 \ (\ap\in\Pi_l); \ (x,\theta)\le 1\} \ ,
\]
and recall that $w=v{\cdot}t_r$, where $v\in W$ and $r\in {Q}^\vee$.

\begin{s}{Proposition} \label{opis-smi} \par
1. $w=v{\cdot}t_r\in \HW^{(s)}_{min}$ \ $\Longleftrightarrow
\left\{ \begin{array}{l} w \text{ is dominant},  \\
                         v(r)\in D^{(s)}_{min}\cap {Q}^\vee \ .
       \end{array}\right.$

2. The mapping $\HW^{(s)}_{min} \to D^{(s)}_{min}\cap Q^\vee$, $w=v{\cdot}t_r\mapsto v(r)$,
is a bijection.
\end{s}\begin{proof}
The argument is the same as in Proposition~\ref{opis-mi}, taking into
account that the constraints for $D^{(s)}_{min}$ are different.
\end{proof}%
In order to compute the number $\#(D^{(s)}_{min}\cap Q^\vee)$,
we perform the following transformation.
Set $\rho^\vee_s=\sum_{\ap_i\in\Pi_s}\varpi^\vee_i$. It is easily seen
that the translation $t_{\rho^\vee_s}$ takes $D^{(s)}_{min}$ to
the dilated closed fundamental alcove
\[
(g+1)\ov{\gA}=\{ x\in V\mid (x,\ap)\ge 0 \ \forall\,\ap\in\Pi;\
   (x,\theta)\le g+1\} .
\]
Here $g=(\theta, \sum_{\ap_i\in\Pi_s}\varpi^\vee_i)$, i.e., it is the sum
of the {\sl short\/} coordinates of $\theta$ (i.e., those corresponding
to the short simple roots).
It is easy to obtain other formulae for $g$. E.g., $g=(\#\Delta_s)/p=
(\rho_s,\theta^\vee)$, where $\rho_s$ is the half-sum of all
positive short roots. If we want to explicitly indicate that $g$
depends on $\Delta$, we write $g_\Delta$.
\\
Although the above translation may not belong to $\HW$, the very existence
of such a transformation and Lemma~1 in \cite{CP2} show that the
following is true

\begin{s}{Lemma}   \label{g+1}
If $g+1$ and the index of connection of
$\Delta$ are relatively prime, then there is an element of\/
$\HW$ that takes $D^{(s)}_{min}$ to $(g+1)\ov{\gA}$.
\end{s}%
The numbers $g$ for all root systems with roots of different lengths
are as follows:

\begin{center}
\begin{tabular}{c||c|c|c|c}
 $\Delta$ & $\GR{C}{p}$ & $\GR{B}{p}$ & $\GR{F}{4}$ & $\GR{G}{2}$ \\ \hline
 $g$ &  $2p-2$ &     2   &  6 &  3
\end{tabular}
\end{center}
It follows that Lemma~\ref{g+1} always applies and hence
$\#(\AN_s)=\#(\HW^{(s)}_{min})=\# ((g+1)\ov{\gA}\cap Q^\vee)$.
In turn, if $g+1$ is relatively prime with the coordinates of
$\theta$, then this number is computed by Eq.~\re{haiman}.
One sees that the condition of relative primeness
does not hold only for $\GR{G}{2}$.
(However, this case can be studied by hand.) Thus, we obtain

\begin{s}{Theorem}   \label{chislo-smi}
Suppose $|\theta|^2/|\theta_s|^2=2$. Then
\[
   \#(\AN_s)= \prod_{i=1}^p \frac{g+e_i+1}{e_i+1}  \ .
\]     \vskip-1ex
\end{s}%
It is easily seen that $\#(\AN_s)=4$ for $\GR{G}{2}$.
\\
Looking at the factors occurring in the formula of Theorem~\ref{chislo-smi},
one may notice that there is a nice formula for $\AN_s$,
which resembles Eq.~\re{chislo-mi}
and also covers the case of $\GR{G}{2}$. Here is it.
\\[.5ex]
Suppose $e_1 < e_2< \ldots < e_p$ and set $n=\#(\Pi_s)$. Then for any
$\Delta$  we have
\begin{equation}  \label{unified}
  \#(\AN_s)= \prod_{i=1}^n \frac{h+e_i+1}{e_i+1}  \ .
\end{equation}
But it is not clear how to prove this {\sl a priori}.
\\
Changing the role of long and short roots in Definition~\ref{def-smi},
one may define $l$-{\it minimal\/} elements, which are in a one-to-one
correspondence with the {\it long\/} antichains. Since the proofs here are
similar, we state only results.
The simplex associated with the $l$-minimal elements is
\[
   D^{(l)}_{min}=\{x\in V \mid (x,\ap)\ge 0 \ (\ap\in\Pi_s); \quad
(x,\ap)\ge 1 \ (\ap\in\Pi_l); \ (x,\theta)\le 2\} \ ,
\]
and the $l$-minimal elements bijectively correspond to the integral points
of $D^{(l)}_{min}$.
The shift in the direction of $\rho^\vee_l=
\sum_{\ap_i\in\Pi_l}\varpi^\vee_i$
takes $D^{(l)}_{min}$ to $(h+1-g)\gA$. Therefore, if $\Delta$ is not
of type $\GR{G}{2}$, then
\begin{equation}  \label{chislo-lmi}
       \#(\AN_l)=\prod_{i=1}^p \frac{h-g+e_i+1}{e_i+1}  \ .
\end{equation}
Since this number can also be computed as $\#(\AN_s)$ for the dual root
system $\Delta^\vee$,
a relation between $g$ for $\Delta$ and $\Delta^\vee$ emerges.
Namely, $g_\Delta +g_{\Delta^\vee}=h$.
\\
Now we turn to considering a `short' analogue of maximal elements.

\begin{rem}{Definition}   \label{def-sma}
$w\in\HW$ is called {\it $s$-maximal\/}, if
\begin{itemize}
\item[\sf (i)] $w$ is dominant;
\item[\sf (ii)]
if $\ap\in\Pi_s$ and $w^{-1}(\ap)=k\delta+\mu$ with
$\mu\in \Delta$, then $k\le 1$;
\item[\sf (iii)]
if $\ap\in\HP_l$ and $w^{-1}(\ap)=k\delta+\mu$ with
$\mu\in \Delta$, then $k\le 0$.
\end{itemize}
\end{rem}%
Using (i), conditions (ii),\,(iii) can be made precise. If $k=0$,
then $\mu\in\Delta^+$; if $k=1$ in (ii), then $\mu\in -\Delta^+$.
\\[.5ex]
We write $\HW^{(s)}_{max}$ for the set of all $s$-maximal
elements. Notice that $\HW^{(s)}_{max}\subset \HW_{max}$.

As in case of maximal elements, we wish to set up a one-to-one
correspondence between
the $s$-maximal elements and a certain subset of $\AN_s$.
In order to distinguish the right subset we need some preparations.
Recall that, although
$\Delta_s$ is not a sub-root system of $\Delta$, it
is a root system in its own right.
Clearly, $\Delta^+_s$ is the set of positive roots for $\Delta_s$.
Let us write $\Pi(\Delta^+_s)$ for the corresponding set of simple
roots. As $\Pi$ itself, it consists of $p$ roots.
Obviously, $\Pi_s\subset \Pi(\Delta^+_s)$. Other roots in $\Pi(\Delta^+_s)$
are in bijection with $\Pi_l$. Each $\beta\in\Pi_l$
is replaced by a short root as follows.
Let $\ap$ be the closest to $\beta$ (in the sense of the Dynkin
diagram) short simple root. The sum of all simple roots in the string
connecting $\ap$ and $\beta$ is a short root, which is a simple root
for $\Delta^+_s$.
\\[.4ex]
{\sf Warning.} Although we often consider antichains lying in
(certain subsets of) $\Delta^+_s$, it is always meant that the ordering
`$\curle$' is inherited  from the whole of $\Delta^+$.

\begin{s}{Proposition}  \label{bij-sma}
The bijection between $\AD_0$ and $\HW_{max}$ described in
Proposition~\ref{bij-ma} gives rise to a bijection between
$\HW^{(s)}_{max}$ and the short antichains lying in
$\Delta^+_s\setminus \Pi(\Delta^+_s)$.
\end{s}%
\begin{proof}
The correspondence described in Proposition~\ref{bij-ma} attaches to
a maximal element $w$ its first layer ideal, $I^w$. But even if $w$ is
$s$-maximal, the generators of $I^w$  may not be short roots.
So, we do not immediately obtain a required short antichain.
To correct this,
we take $I^w\cap\Delta^+_s$. (It is also the set of weights of a
nilpotent $\be$-stable
subspace of
 $\vts$.) The set of generators
(minimal elements) of $I^w\cap\Delta^+_s$ is a short antichain of
$\Delta^+$, which we attach to $w\in \HW^{(s)}_{max}$.

Now, we are to prove that the resulting antichain lies in
$\Delta^+_s\setminus \Pi(\Delta^+_s)$ and that
this correspondence is really a bijection.

Recall from Section~\ref{old} that $\Xi(I)$ is the set of maximal
elements of $\Delta^+\setminus I$ and that in case of maximal elements
$\Xi(I^w)$ is described in Proposition~\ref{gen-ma}.
That description implies that, for $w\in \HW^{(s)}_{max}$, $\Xi(I^w)$
consists of short roots. Since $\Xi(I^w)$ covers all simple roots
(see Remark~\ref{cover}) and consists of short roots, it also covers
all roots from $\Pi(\Delta^+_s)$. (Use the explicit description of
$\Pi(\Delta^+_s$ given above.)  This means that the short antichain
$\Gamma(I^w\cap\Delta^+_s)$ does not contain roots from
$\Pi(\Delta^+_s)$.

\un{Injectivity}. If $w,w'\in \HW^{(s)}_{max}$ are different, then
$\Xi(I^w)\ne \Xi(I^{w'})$. Since these two sets consist of short roots,
we obviously have $I^w\cap\Delta^+_s\ne I^{w'}\cap\Delta^+_s$.

\un{Surjectivity}. If $\Gamma$ is an antichain of $\Delta^+$ lying in
$\Delta^+_s\setminus \Pi(\Delta^+_s)$, then
take all maximal short roots in
$\Delta^+\setminus I\langle\Gamma\rangle$. More precisely, let $\Xi$ be
the set of
short roots $\mu$ such that if $\nu\in\Delta^+_s$ and $\nu\succ\mu$
then there is $\gamma\in\Gamma$ such that $\nu\curge\gamma$.
Then $\Xi$ is a short antichain
that covers all roots in $\Pi(\Delta^+_s)$ and hence the whole of $\Pi$.
In view of Remark~\ref{cover}, $\Xi$
is of the form $\Xi(I^w)$ for some $w\in \HW_{max}$. Finally, since
$\Xi$ consists of short roots, this $w$ is $s$-maximal.
\end{proof}%
The antichains of $\Delta^+$ lying in $\Delta^+_s\setminus \Pi(\Delta^+_s)$
are said to be {\it strictly $s$-positive.}
The corresponding subset of $\AN_s$ is denoted by $\AN_{s,0}$. \\
Once again, the next part of our program is a geometric description.
Set \\[.6ex]
\hbox to \textwidth{ \hfil
  $D^{(s)}_{max}=\{x\in V \mid (x,\ap)\le 1 \ (\ap\in\Pi_s); \quad
(x,\ap)\le 0 \ (\ap\in\Pi_l); \ (x,\theta)\ge 1\}$, \hfil }

\begin{s}{Proposition} \label{opis-sma} \par
1. $w=v{\cdot}t_r\in \HW^{(s)}_{max}$ \ $\Longleftrightarrow
\left\{ \begin{array}{l} w \text{ is dominant},  \\
                         v(r)\in D^{(s)}_{max}\cap {Q}^\vee \ .
       \end{array}\right.$

2. The mapping $\HW^{(s)}_{max} \to D^{(s)}_{max}\cap Q^\vee$,
$w=v{\cdot}t_r\mapsto v(r)$, is a bijection.
\end{s}\begin{proof}
The argument is the same as in Proposition~\ref{opis-mi}, taking into
account that the constraints for $D^{(s)}_{max}$ are different.
\end{proof}%
The translation in the direction of $-\rho^\vee_s$, which belongs to
extended affine Weyl group, takes $D^{(s)}_{max}$ to $(1-g)\gA$.
Since $g-1$ is always relatively prime with the index of connection,
there is also an element of $\HW$ that does the same, cf. Lemma~\ref{g+1}.
As in case of $s$-minimal
elements, we have $g-1$ is relatively prime with the coordinates
of $\theta$, if $\Delta$ is not of type $\GR{G}{2}$.
Therefore, if $\Delta\in\{\GR{B}{p},\,\GR{C}{p},\,\GR{F}{4}\}$, then
\begin{equation}  \label{chislo-sma}
     \#(\AN_{s,0})=\#(\HW^{(s)}_{max})=\prod_{i=1}^p \frac{g+e_i-1}{e_i+1}  \ .
\end{equation}
For $\GR{G}{2}$, this set consists of 2 elements.


\section{Short antichains and the semi-Catalan arrangement}  \label{geners}
\setcounter{equation}{0}

\noindent
In this section, we study a hyperplane arrangement in $V$ that has the same
connection with short antichains in $\Delta^+$ as the Catalan arrangement
has with all antichains. This provides yet another approach
to counting the short and strictly $s$-positive antichains.

\begin{rem}{Definition}  \label{semiCat}
1. The {\it semi-Catalan arrangement\/} in $V$,
$\cats$,
consists of the hyperplanes
$\gH_{\mu,k}$ ($\mu\in\Delta^+_s,\ k=-1,0,1$) and $\gH_{\nu,0}$
($\nu\in\Delta^+_l$).

2. The {\it $m$-semi-Catalan arrangement\/}
in $V$, $\catsm$, consists of the hyperplanes
$\gH_{\mu,k}$ ($\mu\in\Delta^+_s,\ k=-m,\dots,-1,0,1,\dots,m$) and
$\gH_{\nu,0}$ ($\nu\in\Delta^+_l$).
\end{rem}%
All these arrangements are deformations of the Coxeter arrangement.
Notice also that $\cts{0}$ is the usual Coxeter arrangement,
and $\cats=\cts{1}$. \\
First, we are interested in the dominant regions of $\cats$ and their
relation to short antichains. Define a mapping
\[
\psi:\AN_s\to \{\text{the dominant regions of $\cats$}\}
\]
as follows. For $\Gamma\in\AN_s$, let
\[
\Gamma \overset{\psi}{\mapsto} R_\Gamma^{(s)}:=
\{x\in \gC \mid  
(x,\mu)>1, \ \text{ if } \mu\in I\langle\Gamma\rangle\cap\Delta^+_s
\ \& \
(x,\mu)<1,  \ \text{ if } \mu\in \Delta^+_s\setminus I\langle\Gamma\rangle
\} \ .
\]

\begin{s}{Theorem}   \label{bij-reg}
\begin{itemize}
\item[\sf (i)] The mapping $\psi$ is well-defined, and it is a bijection;
\item[\sf (ii)] $R_\Gamma^{(s)}$ is bounded if and only if\
$\Gamma\in \AN_{s,0}$.
\end{itemize}
\end{s}\begin{proof} (i) \
1. Regarding $\Gamma$ as "usual" antichain, we can construct
a region $R_{I\langle\Gamma\rangle}$, as prescribed by Eq.~\re{bij-shi}.
Obviously, $R_{I\langle\Gamma\rangle}\subset R_\Gamma^{(s)}$.
Hence the latter is non-empty.
\\[.4ex]
2. Since the definition of the set $R_\Gamma^{(s)}$ includes a constraint
for
any hyperplane in $\cats$ meeting $\gC$, $R_\Gamma^{(s)}$
cannot contain more than one region.
It is also clear that $R_\Gamma^{(s)}\ne R_{\Gamma'}^{(s)}$,
if $\Gamma\ne\Gamma'$. For, if $\gamma\in \Gamma\setminus
I\langle\Gamma'\rangle$, then $\gH_{\gamma,1}$ separates $R_\Gamma^{(s)}$
and $R_{\Gamma'}^{(s)}$. Hence $\psi$ is injective.
\\[.4ex]
3. The surjectivity of $\psi$ follows from the existence of the inverse
map. Given a region $R$, take the set of walls of $R$ separating
$R$ from the origin. Then the corresponding set of roots
form a short antichain.
\par
(ii) \  If $\Pi(\Delta^+_s)\cap I\langle\Gamma\rangle=\varnothing$, then
$R_{\Gamma}^{(s)}$ belong to the bounded domain
$\{x\in \gC \mid (x,\mu)<1, \ \mu\in \Pi(\Delta^+_s)  \}$.

Conversely, assume $\beta\in I\cap \Pi(\Delta^+_s)$. Recall from
Section~\ref{new} that $\Pi(\Delta^+_s)$ is in bijection with $\Pi$ \
($\beta$ either belong to $\Pi_s$ or is obtained
via a simple procedure from a long simple root). Let $\beta'$
be the simple root in $\Pi$ corresponding to $\beta$ and
$\vp_{\beta'}$ the respective
fundamental weight of $\g$.
Then we claim that if $x\in R_\Gamma^{(s)}$, then $x+t\vp_{\beta'}\in
R_\Gamma^{(s)}$ for any $t\in {\Bbb R}_{\ge 0}$. Indeed, $\beta$ is the minimal
short having nonzero $\beta'$-coordinate. Therefore all short roots
having nonzero $\beta'$-coordinate are in $I\langle\Gamma\rangle$.
This means that $R_\Gamma^{(s)}$ has no upper bound in the direction of
$\vp_{\beta'}$. Thus, $R_\Gamma^{(s)}$ unbounded.
\end{proof}%
Let us look at the relationship between $s$-minimal and $s$-maximal elements
on one hand, and dominant regions of $\cats$ on the other hand.

\begin{s}{Proposition}  \label{near}
\begin{itemize}
\item[\sf (i)] \ Suppose $w\in\HW^{(s)}_{min}$, and let $\Gamma\in\AN_s$
be the corresponding antichain. Then
$w^{-1}\asst\gA \subset R^{(s)}_\Gamma$, and it is the nearest to the
origin alcove in $R^{(s)}_\Gamma$.
\item[\sf (ii)] \ Suppose $w\in\HW^{(s)}_{max}$, and let $\Gamma\in\AN_{s,0}$
be the corresponding antichain. Then
$w^{-1}\asst\gA \subset R^{(s)}_\Gamma$, and it is the most distant
from the origin alcove in $R^{(s)}_\Gamma$.
\end{itemize}
\end{s}\begin{proof}
(i)  It was already observed before that
$w^{-1}\asst\gA \subset R_{I\langle\Gamma\rangle}\subset R_\Gamma^{(s)}$.
Suppose we are inside $w^{-1}\asst\gA$. To get in an alcove that is
closer to the origin, we must cross a wall separating $w^{-1}\asst\gA$
from the origin. These walls correspond to the roots $\ap\in \HP$
such that $w^{-1}(\ap) < 0$. But then $w^{-1}(\ap)=-\delta+\mu$,
where $\mu\in\Delta^+_s$. So that having crossed this wall, we get in
another dominant region of $\cats$.

(ii) Suppose $w\in\HW^{(s)}_{max}$ and we are inside $w^{-1}\asst\gA$.
To get in an alcove that is more distant from the origin, we must cross a
wall of $w^{-1}\asst\gA$ that does not separate $w^{-1}\asst\gA$
from the origin. These walls correspond to the roots $\ap\in \HP$
such that $w^{-1}(\ap) > 0$. In view of the definition of $s$-maximal
elements, there are two possibilities: \
(a) if $w^{-1}(\ap)=\mu\in\Delta+$, then crossing such a wall we
leave $\gC$;
(b) if $w^{-1}(\ap)=\delta-\mu$, $\mu\in\Delta^+_s$, then crossing such
a wall we get in another dominant region of $\cats$.
Hence $w^{-1}\asst\gA$ is the most distant from the origin alcove in
a certain region.
\\
The hyperplanes of $\cats$ separating $w^{-1}\asst\gA$ from the origin
(not necessarily
walls of $w^{-1}\asst\gA$) correspond to the short roots $\mu$ such that
$w(\delta-\mu)<0$, i.e., these roots are exactly the short roots in
the first layer ideal of $w$. According to the correspondence described
in the proof of Proposition~\ref{bij-sma}, the minimal elements of this
set form the short antichain $\Gamma$ attached to $w$.
Thus, $w^{-1}\asst\gA$ lies in the required alcove.
\end{proof}%
Theorem~\ref{bij-reg} implies that the number of short or strictly
$s$-positive antichains can be
found through the use of the characteristic polynomial of
$\cats$. In fact, we are able to compute the characteristic
polynomial for $\catsm$ with any $m$. One should just follow the scheme of
Athanasiadis' proof for the usual $m$-Catalan arrangement, see
\cite[Theorem\,3.1]{ath02}.
Let $P^\vee$ be the coweight lattice and $f=[P^\vee : Q^\vee]$.
(Hence $f$ is the index of connection of $\Delta$.)

\begin{s}{Theorem}    \label{char-pol}
Suppose $t\in \Bbb N$, $t>mg$,
and both $t$, $t-mg$ are relatively prime with all the
coordinates of $\theta$. Then
\[
  \chi(\catsm,t)=\frac{\# W}{f} \#((t-mg)\gA\cap P^\vee) \ .
\] \vskip-1ex
\end{s}\begin{proof}
We give only a sketch of the proof, where we indicate essential
distinction from Athanasiadis' proof for
an $m$-Catalan arrangement, referring to \cite{ath02} for all details.

Let ${\goth P}$ denote the fundamental parallelepiped
$\{ \sum_{i=1}^p y_i\varpi^\vee_i \mid 0\le y_i\le 1\}$.
Set ${\goth P}_t={\goth P}\cap \frac{1}{t}P^\vee$.
Also, let $V^m_{\Delta,t}$ be the set of hyperplanes \\[.5ex]
\centerline{
$\gH_{\mu,k+\frac{n}{t}}$ ($n,k\in\Bbb Z$, $|n|\le m$, $\mu\in\Delta_s$)
\quad and \quad
$\gH_{\gamma,k}$ ($k\in\Bbb Z$, $\gamma\in\Delta_l$). }
\vskip.5ex\noindent
Note that fractional indices are allowed only for hyperplanes orthogonal
to short roots, so that our $V^m_{\Delta,t}$ is different from that of
Athanasiadis.

Given an arrangement $\ca$ in $V$, according to a general result
(Athanasiadis, Bj\"orner-Ekedahl), the
value $\chi(\ca,t)$ is equal to the number of points in
the complement of all hyperplanes, counted after reduction modulo $t$,
i.e.,  in $({\Bbb Z}_t)^p$. More precisely, this equality holds for
infinitely many $t$ (this can be made even more precise, see
\cite[Sect.\,2]{ath02}).
\\[.3ex]
In our situation, this general principle leads to the equality
$\chi(\catsm,t)=\#\{{\goth P}_t\setminus V^m_{\Delta,t}\}$.
Using the fact that ${\goth P}$ contains $\#W/f$ alcoves,
this can be transformed to
\[
  \chi(\catsm,t)=\frac{\#W}{f}\#((\gA\cap \frac{1}{t}P^\vee)\setminus
V^m_{\Delta,t})=
\frac{\#W}{f}\cdot\#((t\gA\cap P^\vee)\setminus
tV^m_{\Delta,t}) \ .
\]
It easily follows from the definition of $V^m_{\Delta,t}$ that
$(t\gA\cap P^\vee)\setminus tV^m_{\Delta,t}$ is obtained from
$t\gA\cap P^\vee$ by deleting the coweights lying on the hyperplanes
$\gH_{\ap,i}$ with $\ap\in\Pi_s$ and $1\le i\le m$.
That is, the set in question is equal to
\[
\{x\in P^\vee \mid (x,\ap)>m \ (\ap\in\Pi_s); \
(x,\ap)>0 \ (\ap\in\Pi_l); \ (x,\theta)<t \} \ .
\]
Finally, the translation by the negative of $m\rho^\vee_s$ (which lies
in $P^\vee$) takes this set to the points of $P^\vee$ lying in the open
simplex $(t-gm)\gA$.
\end{proof}%
Let us discuss consequences of this result.
We use values of $g$ given in Section~\ref{new}.
If $\Delta\in \{\GR{B}{p},\GR{C}{p},\GR{F}{4}\}$ and
$t$ is relatively prime with the coordinates of $\theta$, then
the same holds for $t-mg$ with any $m$. It follows that
\[
   \chi(\catsm,t)=\chi(\cts{0},t-mg)=\prod_{i=1}^p(t-mg-e_i) \ .
\]
(The first equality holds for infinitely many values of $t$; hence it
holds always, as both parts are polynomials in $t$. The second equality is
a statement about Coxeter arrangements.)
In particular,
\begin{equation} \label{cats-bcf}
  \chi(\cats,t)=\prod_{i=1}^p(t-g-e_i) \ .
\end{equation}
Combining Theorems~\ref{zaslavsky} and \ref{bij-reg}, we conclude that
\[
      \#(\AN_s)=|\chi(\cats,-1)|/\#W \quad \text{ and }\quad
      \#(\AN_{s,0})=|\chi(\cats,1)|/\#W \ ,
\]
which coincides, of course, with the formula in Theorem~\ref{chislo-smi}
and Eq.~\re{chislo-sma}.
\\[.5ex]
For $\GR{G}{2}$, we have $g=3$ and
the assumption of Theorem~\ref{char-pol} is satisfied
only if $m$ is even. Therefore
\[
 \chi(\Cat^{2m}_s(\GR{G}{2}),t)=\chi(\Cat^{0}_s(\GR{G}{2}),t-6m)=
(t-6m-1)(t-6m-5) \ .
\]
Using ad hoc arguments, one may derive the ``odd'' formula
\begin{equation} \label{cats-g2}
  \chi(\Cat^{2m+1}_s(\GR{G}{2}),t)=(t-6m-5)(t-6m-7) \ .
\end{equation}
It is also easy to compute $\chi(\Cat^{1}_s(\GR{G}{2}),t)$
directly from the definition of a characteristic polynomial.
\\[.5ex]
Again, it is noteworthy that formulae \re{cats-bcf},\,\re{cats-g2}
for $\chi(\cats,t)$
admit a uniform presentation for all non-simply laced cases,
cf. Eq.~\re{unified}.

\begin{s}{Theorem}  \label{uniform-chi}
If $n=\#(\Pi_s)$ and the exponents of $\Delta$ are increasingly ordered,
then
\[
  \chi(\cats,t)=\prod_{i\le n} (t-h-e_i)\prod_{i\ge n+1}(t-e_i) \ .
\]  \vskip-1ex
\end{s}%
Of course, it would be interesting to have a uniform proof
(explanation) for this.

\begin{rem}{Remarks}   \label{semi-shi}
1. One may consider `short' analogues for other arrangements
associated with root systems. For instance,
the {\it extended semi-Shi arrangement\/},
$\text{Shi}_s^m(\Delta)$,
is the collection of hyperplanes
$\gH_{\mu,k}$ ($\mu\in\Delta^+_s,\ k=-m+1,\dots,-1,0,\dots,m$)
and $\gH_{\nu,0}$ ($\nu\in\Delta^+_l$).
It is not hard to compute that, for $\GR{C}{2}$,
the characteristic polynomial is equal to $(t-2m-1)^2$.
For $\GR{G}{2}$, it is equal to $(t-3m-1)(t-3m-2)$, at least
if $m\ge 3$. I conjecture that the following formula holds in
general:
\[
  \chi(\text{Shi}_s^m(\Delta),t)=\prod_{i=1}^p
 (t-mg_{\Delta^\vee}-e_i(\Delta_l)) \ ,
\]
where $\{e_i(\Delta_l)\}$ are the exponents for the root system $\Delta_l$.
For instance, in case of $\GR{F}{4}$ we have $g_\Delta=g_{\Delta^\vee}=6$
and $\Delta_l$ is of type $\GR{D}{4}$. Therefore the conjectural expression
is $(t-6m-1)(t-6m-3)^2(t-6m-5)$.

2. The dominant regions of $\cats$ provide a connection, in the spirit
of \cite{gs}, with nilpotent $G$-orbits in $\vts$. I hope to
discuss this topic in a forthcoming publication.
\vskip-1ex
\end{rem}


\section{Some numerical complements}  \label{numer}
\setcounter{equation}{0}

\noindent
In this section, we collect several results that can be proved in a
case-by-case fashion.

\subsection{}
We know the number of all and short antichains for all irreducible
reduced root systems. Using this, on may observe that $\#(\AN_s)$
divides $\#(\AN)$
in all cases. Furthermore, the ratio has, {\it a posteriori\/}, an interesting
description. Namely, let $\Delta(\Pi_l)$ be the root system whose
set of simple roots is $\Pi_l$. Notice that $\Delta(\Pi_l)$ is smaller
than $\Delta_l$, and that the former is irreducible, since
$\Pi_l$ is a connected subset of the Dynkin diagram.
Write $\AN(\Pi_l)$ for the set of all antichains in
$\Delta(\Pi_l)^+$.

\begin{s}{Theorem}  \label{product}
 $\#(\AN)=\#(\AN_s){\cdot} \#(\AN(\Pi_l))$.
\end{s}\begin{proof}
Case-by-case verification. For instance, in case of $\GR{F}{4}$
we have $\#(\AN)=105$, $\#(\AN_s)=21$, and $\Delta(\Pi_l)$ is of
type $\GR{A}{2}$, where one has five antichains.
\end{proof}%
Of course, this proof is not illuminating. One may consider
a natural mapping $\AN \to \AN_s$ that takes
$\Gamma$ to the set of minimal elements of $I\langle\Gamma\rangle
\cap \Delta^+_s$. For $\GR{C}{p}$,
all fibres of this mapping has the same cardinality,
which is 2.
To some extent, this is an explanation in this case.
Unfortunately, the "equicardinality" property does not hold
for $\GR{F}{4}$ and $\GR{G}{2}$.
The statement of Theorem~\ref{product} can be compared with another equality,
which is easy to prove. The reflection $s_\gamma\in W$ is called short, if
$\gamma\in\Delta^+_s$. Let $W_s$ be the (normal) subgroup of $W$ generated
by all short reflections, and let $W(\Pi_l)$ be the Weyl group of
$\Delta(\Pi_l)$. Then $W\simeq W_s\rtimes W(\Pi_l)$ \ (a semidirect
product).

\subsection{}
We have shown that the short antichains of $\Delta^+$ lying in
$\Delta^+_s \setminus\Pi(\Delta^+_s)$
are in a one-to-one correspondence with $s$-maximal elements,
and then computed their number. However, it is also natural to
enumerate the short antichains lying in
$\Delta^+_s \setminus\Pi_s$. (Recall that $\Pi_s$ is a proper
subset of $\Pi(\Delta^+_s)$.) Set
$\AN_{ss}=\{\Gamma\in\AN_s\mid \Gamma\cap \Pi_s=\varnothing\}$.
I did not find a suitable bijection
for $\AN_{ss}$, but the following formula for the cardinality
is true:
\begin{equation}  \label{bezpi_s}
 \#(\AN_{ss})=
  \prod_{i=1}^n\frac{h+e_i-1}{e_i+1} \ ,
\end{equation}
where the notation is the same as in Theorem~\ref{uniform-chi}.
Again, this formula bears a striking resemblance with
Theorem~\ref{chislo-ma}. Direct calculations show that this gives us
the correct
number for $\GR{B}{p}$ (this is easy, because there is only a few short
roots), $\GR{F}{4}$,\,$\GR{G}{2}$.
\\[.4ex]
The argument for $\GR{C}{p}$ goes as follows.
The set of positive roots $\Delta^+(\GR{C}{p})$ is naturally
represented by the shifted Ferrers diagram of shape
$(2p-1,2p-3,\dots,1)$, and the ideals
are represented by suitable subdiagrams of it,
see slightly different versions in
\cite[Sect.\,2]{shi},\,\cite[Sect.\,3]{CP1},\,\cite[Sect.\,5]{duality}.
In these interpretations, the long roots are represented by the boxes
in an extreme diagonal of this shifted Ferrers diagram, and the simple
roots correspond to the boxes of another ("opposite") diagonal.
These two diagonals have a unique common box, corresponding to the long
simple root. If we want
to obtain an ideal whose generators are short and
contain no short simple roots, then
we just erase both these diagonals and consider a subdiagram of
the smaller shifted diagram. But this smaller shifted Ferrers diagram,
which is of shape $(2p-3,2p-5,\dots,1)$,
can be thought of as the set of positive roots for $\GR{C}{p-1}$.
Thus, the  number $\#(\AN_{ss})$ for $\GR{C}{p}$
equals the number $\#(\AN)$ for $\GR{C}{p-1}$.
The latter is known to equal $\genfrac{(}{)}{0pt}{}{2p-2}{p-1}$, which
is consistent with Eq.~\re{bezpi_s}. Actually, we obtain more.
Our bijection between
$\AN_{ss}(\GR{C}{p})$ and $\AN(\GR{C}{p-1})$ preserves the number of
elements. Therefore, we conclude that
the number of $k$-element antichains  in $\AN_{ss}(\GR{C}{p})$ is equal to
$\genfrac{(}{)}{0pt}{}{p-1}{k}^2$, $k=0,1,\dots,p-1$.

\subsection{}
Counting antichains with respect to the number of generators
yields an interesting $q$-analogue of $\#(\AN)$, see
\cite{ath03},\,\cite{duality}. In case of two root lengths, one may
consider, of course, a 2-parameter refinement.
Set
\[
  \AN\langle k,m\rangle=\{\Gamma\in\AN\mid
\#(\Gamma\cap\Delta^+_s)=k \ \ \& \ \
\#(\Gamma\cap\Delta^+_l)=m\},\quad  a_{k,m}=\#\AN\langle k,m\rangle \ ,
\]
and consider the generating function $\cF(t,u)=\sum_{k,m\ge 0}a_{k,m}t^ku^m$.
We have
\begin{itemize}
\item[$\GR{G}{2}$:] \  $\cF(t,u)=1+3t+3u+tu$.
\item[$\GR{F}{4}$:] \  $\cF(t,u)=1+12t+12u+8t^2+39tu+8u^2+12t^2u+
12tu^2+t^2u^2$.
\end{itemize}
The symmetry of these polynomials stems from the fact the corresponding
root systems are self-dual. Since the root systems of type {\bf B} and
{\bf C}
are dual to each other, the corresponding matrices $(a_{k,m})$
are mutually transposed. So, it suffices to handle the case of
$\GR{C}{p}$. Each pair of long roots in $\Delta^+(\GR{C}{p})$ is
comparable, hence any antichain contains at most one long root.
So that we are to determine the coefficients
$a_{k,0},a_{k,1},\ (k=0,1,\dots,p-1)$. In \cite[Sect.\,5]{duality},
we constructed an involution on the set $\AN(\GR{C}{p})$, which maps
$\AN\langle k,0\rangle$ onto $\AN\langle p-1-k,1\rangle$.
Hence $a_{p-k-1,1}=a_{k,0}$ and we have to only  count the number of
short antichains with $k$ elements. Using shifted Ferrers diagrams,
it can be shown that
$a_{k,0}=
\genfrac{(}{)}{0pt}{}{p}{k}\genfrac{(}{)}{0pt}{}{p-1}{k}$.
(In this situation, short simple roots are allowed, so that one has to
erase only one diagonal and work with the shifted Ferrers diagram
of shape $(2p-2,2p-4,\dots,2)$.)

\end{document}